%% file: aliasing_main.tex
\documentclass[preprint,12pt]{elsarticle}
\pdfoutput=1
\usepackage{graphicx}
\usepackage{amssymb}
\usepackage{amsthm}
\usepackage{graphicx}
\usepackage{graphics}
\usepackage{upgreek}
\usepackage{float}
\usepackage{epstopdf}
\usepackage{amsmath}
\usepackage{wrapfig}
\usepackage{pdflscape}
\usepackage{listings}
\usepackage{subcaption}
\usepackage{multirow}
\usepackage[percent]{overpic}
\usepackage{multicol}
\usepackage{color}
\usepackage{stmaryrd}
\usepackage[hidelinks]{hyperref}

\usepackage{diagbox}

\usepackage{booktabs}
\usepackage{epstopdf}
\usepackage[section]{placeins} % For '\input{file}'

\usepackage{calc}
\newlength\myheight
\newlength\mydepth
\settototalheight\myheight{Xygp}
\newcommand*\inlinegraphics[1]{%
  \settototalheight\myheight{Xygp}%
  \settodepth\mydepth{Xygp}%
  \raisebox{-\mydepth}{\includegraphics[height=\myheight]{#1}}%
}
\newcommand\orcid[1]{\href{https://orcid.org/#1}{\inlinegraphics{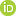}}}

\usepackage[titletoc]{appendix}

\usepackage[section]{algorithm}
\makeatletter
\def\BState{\State\hskip-\ALG@thistlm}
\makeatother

\journal{}

\graphicspath{ {Figs/} }

\newcommand{\etal}{et~al.}

\newcommand{\half}{\frac{1}{2}}
\newcommand{\bigo}{\mathcal{O}}
\newcommand\bigob[1]{\mathcal{O}(#1)}

\newcommand{\hfd}{\hat{f}^{\delta}}

\newcommand{\hfI}{\hat{f}^{\delta I}}
\newcommand{\hfC}{\hat{f}^{\delta C}}
\newcommand{\hfD}{\hat{f}^{\delta D}}

\newcommand{\hud}{\hat{u}^{\delta}}

\newcommand\hb[1]{\hat{\mathbf{#1}}}

\newcommand\px[2]{\frac{\partial #1}{\partial {#2}}}

\newcommand\dx[2]{\frac{\mathrm{d} #1}{\mathrm{d} #2}}
\newcommand\dxi[3]{\frac{\mathrm{d}^{#1}#2}{\mathrm{d} {#3}^{#1}}}
\newcommand\rint[2]{\int^{1}_{-1}{#1} \mathrm{d}{#2}}

\newcommand\CONDITION[2]%
\begin{document}

\begin{frontmatter}

\title{All Equalities Are Equal, but Some Are More Equal Than Others:\\ 
\large{The Effect of Implementation Aliasing on the Numerical Solution to Conservation Equations}}

\author{W. Trojak\corref{cor1}~\orcid{0000-0002-4407-8956}}
\ead{wt247@cam.ac.uk}
\cortext[cor1]{Corresponding author}
\address{Department of Engineering, University of Cambridge, CB2 1PZ}

\author{A. Scillitoe\corref{cor2}~\orcid{0000-0001-8971-7224}}
\ead{ascillitoe@turing.ac.uk}
\address{Alan Turing Institute, Kings Cross, London, UK, NW1 2DB}

\author{R. Watson\corref{cor2}~\orcid{0000-0001-9683-8489}}
\ead{r.watson@qub.ac.uk}
\address{School of Mechanical and Aerospace Engineering, Queen's University Belfast, Belfast, UK, BT9 5AH}

\begin{abstract}
	 We investigate the effect of aliasing when applied to the storage of variables, and their reconstruction for the solution of conservation equations. In particular, we investigate the effect on the error of storing primitives versus conserved variables for the Navier-Stokes equations. It was found that storing the conserved variables introduces less dissipation and that the dissipation caused by constructing the conversed variable from the primitives grows factorially with the order. Hence, this problem becomes increasingly important with the continuing move towards higher orders. Furthermore, the method of gradient calculation is investigated, as applied to the viscous fluxes in the Navier-Stokes equations. It was found that in most cases the difference was small, and that the product rule applied to the gradients of the conserved variables should be used due to a lower operation count. Finally, working precision is investigated and found to have a minimal impact on free-stream-turbulence-like flows when the compressible equations are solved, except at low Mach numbers.
\end{abstract}


\end{frontmatter}

%% Start line numbering here if you want
%\linenumbers

% ================================================================================

%% main text
\input{introduction}
\input{fr}
\input{palias}
\input{variables}
\input{icv}
\input{tgv}

\section{Conclusions}\label{sec:conclusions}
	The effect of implementation on the accuracy of solving conservative PDEs was investigated. It was found that the aliasing error when constructing conserved variables from primitive variables is factorially dependant on the order of the method. This result was confirmed via numerical tests, where the implementation error becomes more apparent as the order was increased. It was also found that storing the primitive variables and constructing the conserved variables led to higher dissipation, and therefore it is recommended that the conserved variables are stored, as this method is also faster.
	
	Two methods were also tested for calculating the gradients required when forming the viscous fluxes for the Navier-Stokes equation, and it was found that the different methods were similar due to the small absolute value of the viscous flux. Therefore for most applications using the product rule on the gradient of the conserved variables is recommended, as it reduces the operation count. Finally, the effect of working precision was investigated, and it was found that single precision gives acceptable results in pseudo-free-stream-turbulence flows. The difference caused by precision was found to be more pronounced in low Mach number regimes for compressible schemes.

\section*{Acknowledgements}
\label{sec:ack}
	The support of the Engineering and Physical Sciences Research Council of the United Kingdom is gratefully acknowledged under the award reference 1750012. The authors would also like to thank Nvidia for the GPU Seeding Grant received.

\section*{References}
\bibliographystyle{elsarticle-num}
\bibliography{library}

%% Authors are advised to submit their bibtex database files. They are
%% requested to list a bibtex style file in the manuscript if they do
%% not want to use model1-num-names.bst.

\clearpage
\begin{appendices}
\input{nomenclature}
\end{appendices}

% Show the list of todos in the document.  Needed to avoid stupid warnings/errors when using the todo package
%\todos

\end{document}

%% file: introduction.tex
%!TEX root = ./aliasing_main.tex
\section{Introduction}
\label{sec:intro}	

	Over the course of the past three decades, Large Eddy Simulation (LES) has become increasingly used for the exploration of flow physics. Looking forward to how CFD will be used tomorrow, NASA's CFD Vision 2030~\cite{Slotnick2014} predicts that hybrid RANS/LES and wall-modelled LES will become increasingly used in aerospace design, and that these methods are likely to prevail until sufficient technological developments allow for wall-resolved LES to become a feasible part of the design process. The effect of this continuing shift from low fidelity modelling to high fidelity simulation is that the gap is bridged, in part, by adapting existing RANS tools for LES. For example, ANSYS Fluent began its life with a RANS turbulence modelling approach~\cite{Wahab2014}.
    
    Early in the development of numerical methods for computationally approximating solutions to PDEs, high order methods became of interest --- as they offered potentially lower mesh requirements and reduced error. For example, the development of Discontinuous-Galerkin method by Reed and Hill~\cite{Reed1973} began a long journey which has proved fruitful due to its super-convergence~\cite{Cockburn1999}. However, with the advent of LES, these high-order methods have become the focus of significant renewed research effort. Multiple tools capable of large high fidelity calculations are becoming available, for example Nektar++\cite{Karniadakis2013}, PyFR~\cite{Witherden2014}, etc. 
    
    Throughout this process of adapting and developing methods, there seems to have been insufficient public consideration paid to the effect implementation details can have on the solution, especially as the methods become more sensitive via high order in the pursuit of higher fidelity. In particular, the investigation presented here is concerned with the form in which the variables are stored when solving a conservation equation. In the context of fluid mechanics, this commonly comes down to the question of whether the primitive or conservative variables are stored. To the knowledge of the authors, the only justification given for this choice one way or the other comes from Fluent's documentation~\cite{Fluent_25_5_2_2006}, where the given arguments are `it is a natural choice when solving incompressible flows' and `to obtain more accurate velocity and temperature gradients in viscous fluxes, and pressure gradients in inviscid fluxes'.
    
    We, therefore, propose investigating this further, by attempting to answer the question: is the aliasing error introduced through the construction of the terms required in fluids dynamics sufficient for one method to be favourable? Furthermore, we propose to explore the second point raised above: is there a significant difference in constructing the gradients required for viscous fluxes when the variables are stored differently?
    
    A further point which we shall briefly explore is, whether the received wisdom that variables should be stored at double precision or whether single precision is, in fact, sufficient. In most typical calculations the norm is to use double precision throughout, however, as the size of problems to be tackled grows, so does the memory usage. It would be beneficial to both reducing memory overhead and increasing computational speed if single precision were used. Further to this, some hardware --- notably, a large number of GPUs --- include only a small number of double precision arithmetic units and therefore the increase in computational speed can be as much as by a factor of 32 as a result of the move from 64 to 32 bit precision. Some investigation into this question has been performed, notably by Homann~\etal~\cite{Homann2007} on the DNS of incompressible homogeneous turbulence using a variable precision incompressible pseudo-spectral scheme. However they saw little to no difference when the precision was changed, but this may have been due to the explicit enforcement of incompressibility. Another investigation into precision was presented in the review paper by Bailey~\cite{Bailey2005}, that spanned several physics regimes. This investigation, however, was in the opposite direction, looking at the effect of 128-bit precision. It was found that it could be important and concluded that better support of adaptive precision should be made by software and hardware. Therefore, we wish to investigate if the same insensitivity to working precision is true for high order polynomial based methods, or like Bailey~\cite{Bailey2005} it could, in fact, be important.
    
    The remainder of this paper is structured around a high-order numerical method introduced in section~\ref{sec:fr}. The theory of polynomial aliasing and the effect of the order is presented in section~\ref{sec:palias}. We then go on to set out in section~\ref{sec:variables} the variable forms and conversion methods that will be investigated. Then, in sections~\ref{sec:icv}~\&~\ref{sec:tgv}, numerical experiments with Euler's equations and the Navier-Stokes equations are performed respectively. This is followed by the variation of the working precision when applied to the Navier-Stokes equations. Finally, the conclusions are presented in section~\ref{sec:conclusions}.

%% file: fr.tex
%!TEX root = ./aliasing_main.tex
\section{High-Order Flux Reconstruction}\label{sec:fr}
	To provide a flexible framework for performing simulation at various orders of accuracy, we will make use of the high-order method, Flux Reconstruction~(FR)~\cite{Huynh2007,Vincent2010}. This section aims to introduce the methodology behind FR, helping to inform the later investigation into aliasing. FR is broadly based on the techniques used in Nodal Discontinuous-Galerkin~\cite{Hesthaven2008}, as such, we begin by subdividing the domain $\mathbf{\Omega}$ into $n$ sub-domains.
	\begin{equation}
		\mathbf{\Omega} = \bigcup^N_{n=1}\mathbf{\Omega}_n, \quad \mathrm{and} \quad \mathbf{\Omega}_i \cap \mathbf{\Omega}_j = \emptyset \: \forall \: i\ne j
	\end{equation}
	If we then focus on the method as applied to 1D conservation equations, we can define a spatial transformation from the physical sub-domain $\mathbf{\Omega}_n \in [x_n,x_{n+1}]$ to a reference domain $\hat{\mathbf{\Omega}} \in [-1,1]$. This can be achieved via the mapping $\Gamma_n: x\rightarrow \xi$, where $x$ is a variable in $\mathbf{\Omega}_n$ and $\xi$ is in $\hat{\mathbf{\Omega}}$. $\Gamma_n$ is then defined as:
	\begin{equation}
		\xi = \Gamma_n(x) = 2\bigg(\frac{x-x_n}{x_{n+1}-x_{n}}\bigg) - 1
	\end{equation}
	If we proceed to solve the 1D first order conservation equation, then:
	\begin{equation}
		\px{u}{t} + \px{f}{x} = 0
	\end{equation}
	where $u$ is the conserved variable and the flux is $f=f(u)$. Within each sub-domain, we use the data stored at a series of points to form a local polynomial of $u$ and $f$:
	\begin{align}
		\hud(\xi) &= \sum^p_{i=0}\hud_il_i(\xi) \label{eq:ud}\\
		\hfD(\xi) &= \sum^p_{i=0}\hfd_il_i(\xi) \label{eq:fd}
	\end{align}
	where $p$ is the order, the superscript delta symbolises that the polynomials are local to one element, and the hat marks that the variable has been transformed from the physical to reference domain. In the case of the flux, there is also a $D$ to symbolise it is currently only a fit based on the data and hence not strictly continuous. Here, the polynomial basis, $l_i(\xi)$, is the Lagrange basis, defined as:
	\begin{equation}
		l_i(\xi) = \prod^{p}_{\substack{j=0 \\ j\ne i}}\frac{\xi-\xi_i}{\xi_j-\xi_i}.
	\end{equation}
	Now the approximation in Eq.(\ref{eq:ud}) can be used to extrapolate to the edges of the element at $\xi=\pm 1$. This data can be combined with the edge values of the surrounding elements to calculate a common value at each element interface. This is key in enabling the solution between elements to be made continuous. There are several methods of finding a common value, for example, central differencing can be used, but at the expense of needing smoothing. Alternatively, a method such as a Riemann solver can be used that accounts for the upwind direction in hyperbolic equations~\cite{Toro2009} and consequently adds some stabilising dissipation. 
	
	With left and right common interface values calculated, defined as $\hfI_L$ and $\hfI_R$, the common value then needs to be propagated into the element to form a continuous solution. This is achieved via correction functions, $h_L$ and $h_R$, that have the following properties:
	\begin{align}
		h_L(-1) = h_R(1) = 1& \\
		h_R(-1) = h_L(1) = 0&.
	\end{align}
	There are several families of correction function, with it being known that the choice can have a large impact on the behaviour of the method~\cite{Vincent2011,Trojak2018,Trojak2018b,Trojak2018f}. The correction to the flux term is calculated as:
	\begin{equation}
		\hfC = (\hfI_L-\hfD_L)h_L + (\hfI_R-\hfD_R)h_R
	\end{equation}
	then formulating the corrected flux gradient:
	\begin{equation}\label{eq:fcorrected}
		\px{\hfd}{\xi} = \px{\hfD}{\xi} + (\hfI_L-\hfD_L)\dx{h_L}{\xi} + (\hfI_R-\hfD_R)\dx{h_R}{\xi}
	\end{equation}
	lastly by using the transformed equation:
	\begin{equation}
		\px{\hud}{t} + \px{\hfd}{\xi} = 0
	\end{equation}
	with Eq.(\ref{eq:fcorrected}), a suitable temporal integration method may be used to advance the solution in time.     

    One advantage of FR is that by varying the number of points and hence the order of the interpolation in Eq.(\ref{eq:ud}~\&~\ref{eq:fd}), the order accuracy of the scheme may be changed with relative ease. Hence, this will easily allow for the effects of aliasing under investigation here to monitored at different orders, with the aim of drawing further conclusions for the development of high-order methods.

%% file: palias.tex
%!TEX root = ./aliasing_main.tex
\section{Discrete Polynomial Aliasing}
\label{sec:palias}

	In this section, we wish to introduce polynomial aliasing and how this error enters the solution of conservative equations. Let us start by studying a simple generalised conservative equation:
	\begin{equation}\label{eq:gen_conv}
		\px{u}{t} + \px{f}{x} = 0
	\end{equation}
	In this case, we will solve on the periodic domain $[-1,1]$, for simplicity. Then the effect of solving this numerically is that we have some finite basis. Let us then say that the solution, $u$, may be constructed as some $p^{\mathrm{th}}$ order polynomial. In this case, the Legendre basis will be used:
	\begin{equation}\label{eq:uleg}
		u = \sum^p_{i=0}\tilde{u}_i\psi_i(x)
	\end{equation}	
	where $\psi_m$ is the $m^{\mathrm{th}}$ order Legendre polynomial of the first kind. If the flux function is then $f=f(u^n)$ for $n\in\mathbb{N}$, then for the flux we get:
	\begin{equation}\label{eq:fleg}
		f = \sum^{np}_{i=0}\tilde{f}_i\psi_i(x)
	\end{equation}
	However, as previously stated, the functional space of the numerical solver is limited to be $p^{\mathrm{th}}$ order. To understand the error let us discretise Eq.(\ref{eq:gen_conv}) by using Eq.(\ref{eq:uleg}~\&~\ref{eq:fleg}). It is then possible to use a key result about Legendre polynomials, that simple truncation of the series gives the least squares projection of the polynomial. Therefore we may write: 
	\begin{equation}
		\px{}{t}\sum^p_{i=0}\tilde{u}_i\psi_i(x) = -\px{}{x}\sum^{p}_{i=0}\tilde{f}_i\psi_i(x)
	\end{equation} 
	The aliasing error, $e_a$, is defined as:
	\begin{equation}
		\px{}{x}\sum^{np}_{i=0}\tilde{f}_i\psi_i = \px{}{x}\sum^{p}_{i=0}\tilde{f}_i\psi_i + e_a
	\end{equation}
	and then going on to define the differentiated polynomial coefficients as:
	\begin{equation}
		\px{}{x}\sum^{p}_{i=0}\tilde{f}_i\psi_i = \px{}{x}\sum^{p-1}_{i=0}\tilde{f}^\prime_i\psi_i
	\end{equation}
	Therefore:
	\begin{equation}
		\sum^{p-1}_{i=0}\tilde{f}^\prime_i\psi_i(x) - \sum^{np-1}_{i=0}\tilde{f}^\prime\psi_i(x) = -\sum^{np-1}_{i=p}\tilde{f}^\prime\psi_i(x) = -e_a
	\end{equation}
	The exact value of the coefficients $\tilde{f}$ and $\tilde{f}^\prime$ are dependent on the numerical method and on $u$, but we may write the energy in the error term as:
	\begin{equation}
		\|e_a\|^2_{L_2} = \rint{\Bigg(\sum^{np-1}_{n=p}\tilde{f}^{\prime}_n\psi_n(x)\Bigg)^2}{x} = \sum^{np-1}_{n=p}\frac{2\big(\tilde{f}^\prime_n\big)^2}{2n+1} \leqslant \frac{2}{2p+1}\|\tilde{f}^\prime\|_{\infty}^2
	\end{equation}
	The aim of this has been to explain that in the context of a numerical method, the aliasing arises due to the inability of the method to resolve higher order terms, and, furthermore, the energy is removed as dissipation due to the sign of $e_a$.

    This example demonstrates this, but it does not demonstrate the exact behaviour this paper is to investigate, that aliasing behaviour depends on the way in which the variables are stored and used. In order to investigate this, we will modify Burgers' equation. First, we will introduce a key result in aliasing. Defining the interpolation remainder:
	\begin{equation}
		\mathcal{R}_pf = f - \mathcal{L}_pf
	\end{equation}
	where $\mathcal{L}_p$ is a $p^\mathrm{th}$ order linear \emph{interpolation} operator. From Kress~\cite{Kress1998} it can then be stated that:
	\begin{equation}
		\mathcal{R}_pf(\xi) = \frac{f^{(p+1)}(\epsilon)}{(p+1)!}\prod^{p}_{i=0}(\xi-\xi_i)
	\end{equation}
	where $\xi_i$ are the interpolation points and $\epsilon$ is dependant on $\xi$. If we take Burgers' equation and square the conserved variable, we get the following:
	\begin{equation}
		\px{u^2}{t} + \px{u^4}{\xi} = 0
	\end{equation}
	This produces an opportunity for information to be stored in two ways, which are comparable to methods used for Euler's equations --- namely, storing $u$ or $u^2$, and forming the flux term by either squaring $u^2$, or by raising $u$ to the power of four. This gives two possible flux polynomials when transformed into the computational domain:
	\begin{align}\label{eq:f2}
		\hat{u}^2(\xi) = \hat{f}_2 &= \sum^{2p}_{i=0}\tilde{f}_{2,i}\psi_i \\
		\hat{u}^4(\xi) = \hat{f}_4 &= \sum^{4p}_{i=0}\tilde{f}_{4,i}\psi_i \label{eq:f4}
	\end{align}
	To understand how errors may then enter the solution, we wish to understand the scaling of the remainder of the flux interpolation to a finite polynomial space of order $p$. The maximal norm can then be used to give an estimate as:
	\begin{equation}
		\|\mathcal{R}_pf\|_{\infty} \leqslant \frac{1}{(p+1)!}\|q_{p+1}\|_{\infty}\|f^{(p+1)}\|_{\infty}.
	\end{equation}
	Here we define $q_{p+1}$ as:
	\begin{equation}
		q_{p+1} = (\xi-\xi_0)(\xi-\xi_1)\dots(\xi-\xi_p)
	\end{equation}
	with $\xi_i$ being the points at which the value of $f$ is stored. Taking the domain to be $[-1,1]$ therefore $\|q_{p+1}\|_{\infty}\leqslant 2$. We now use Eq.(\ref{eq:f2}~\&~\ref{eq:f4}) to refine the remainder estimates, which requires a bounding value of $\|f^{(p+1)}\|_\infty$. Firstly, it is known that the maximum absolute value of a Legendre polynomial is at $\xi=\pm 1$ and, due to the recursive definition of Legendre polynomials, the maximum value of the derivative is at $\xi=\pm 1$. If the value of a differentiated Legendre polynomial at $\pm1$ is:
	\begin{equation}\label{eq:legd_edge}
		\dxi{m}{\psi_n(\pm1)}{\xi} = \frac{(\pm1)^{n-m}(n+m)!}{2^nn!(n-m)!}
	\end{equation}
	A consequence is that, for a given set of differentiated Legendre polynomials, $\{\psi_n^\prime,\psi_n^{\prime\prime},\dots,\psi_n^{(m)}\}$, the maximum value in this set is the edge value of the $m^\mathrm{m}$ derivative. Hence a bound can be placed on $\|f^{(p+1)}\|_\infty$ using Eq.(\ref{eq:legd_edge}) and the maximum Legendre mode coefficient as:
	\begin{align}
		\|f_2^{(p+1)}\|_\infty &\leqslant \bigg[\frac{2(3p+1)!}{2^{2p}(p)!(2p)!}\bigg]\max_{i\in\{0\dots2p\}}{|\tilde{f}_{2,i}|} \\
		\|f_4^{(p+1)}\|_\infty &\leqslant \bigg[\frac{2(5p+1)!}{2^{4p}(3p)!(4p)!}\bigg]\max_{i\in\{0\dots4p\}}{|\tilde{f}_{4,i}|}
	\end{align}
	Hence, the interpolation remainder may be bounded as:
	\begin{align}\label{eq:r2}
		\|\mathcal{R}_pf_2\|_\infty &\leqslant 4\bigg[\frac{(3p+1)!}{2^{2p}(p)!(2p)!(p+1)!}\bigg]\max_{i\in\{0\dots2p\}}{|\tilde{f}_{2,i}|} \\
		\|\mathcal{R}_pf_4\|_\infty &\leqslant 4\bigg[\frac{(5p+1)!}{2^{4p}(3p)!(4p)!(p+1)!}\bigg]\max_{i\in\{0\dots4p\}}{|\tilde{f}_{4,i}|}\label{eq:r4}
	\end{align}
	It can be proved by induction that:
	\begin{equation}
		\frac{(3p+1)!}{2^{2p}(p)!(2p)!} \leqslant \frac{(5p+1)!}{2^{4p}(3p)!(4p)!}, \quad \forall p\in \mathbb{N}
	\end{equation}
	From this, there are two conclusions that can be drawn. Firstly, the interpolation remainder of $f_4$ will always be bigger than $f_2$. Secondly, the difference between the remainders will grow factorially fast as the order is increased. Therefore, higher order methods will be greatly more affected by this mechanism of error introduction. 
	
	Now considering FR in 1D for a conservative equation we get:
	\begin{align}
		\px{\hud}{t} + \px{\hfd}{\xi} &= 0 \\
		\px{\hfd}{\xi} &= \px{\hfD}{\xi} + (\hfI_L - \hfD_L)\dx{h_L}{\xi} + (\hfI_R - \hfD_R)\dx{h_R}{\xi}.
	\end{align}
	Therefore, we are concerned with two extensions of the remainder derived earlier: the interpolated remainder of the interpolation, and the error of the interpolation to the left and right interfaces. To calculate this let us write:
	\begin{equation}
		\mathcal{R}^{\prime}_pf = \dx{f}{x} - \dx{}{x}\mathcal{L}_pf = \dx{}{x}\mathcal{R}_pf
	\end{equation}
	Hence we can differentiate the result of Kress~\cite{Kress1998}:
	\begin{equation}
		\mathcal{R}^{\prime}_pf(\xi) = \frac{f^{(p+1)}(\epsilon)}{(p+1)!}\dx{}{\xi}\prod^p_{i=0}(\xi-\xi_i), \quad \xi\in[-1,1]
	\end{equation}
	where $\epsilon \in [-1,1]$. Finally we can write:
	\begin{equation}
		\|\mathcal{R}^{\prime}_pf\|_{\infty} = \frac{p}{(p+1)!}\max_{i\in\{0..p\}}(\|q_{p}^i\|_{\infty})\|f^{(p+1)}\|_{\infty} 
	\end{equation}
	defining $q^i_p$ as:
	\begin{equation}
		q^i_p = \prod^{p}_{j=0,j\ne i}(\xi-\xi_j)
	\end{equation}
	If we then apply the results of Eq.(\ref{eq:r2}~\&~\ref{eq:r4}), it is demonstrated that the primary difference in the gradient remainder is a factor of $p$. If we then consider the more straightforward case of the remainder from interface interpolation, \emph{i.e.} $\mathcal{R}_pf(\pm1)$. 
    \begin{equation}\label{eq:rlr}
		\mathcal{R}_pf(\pm1) = \frac{f^{(p+1)}(\epsilon)}{(p+1)!}\prod^p_{i=0}(\pm1-\xi_i)
	\end{equation}

	We will not consider the infinity norm in this case, as Eq.(\ref{eq:rlr}) gives sufficient details. The primary feature of note is that the behaviour of this remainder is primarily influenced by the interpolation point locations. For example, if $\xi_0 =-1$ and $\xi_p=1$, as in a Gauss-Lobatto quadrature the aliasing error introduced through this mechanism would be zero. However, for other reasons explored by Castonguay~\cite{Castonguay2011} this is problematic in higher dimensions. By comparing the remainder due to differentiation and interface interpolation, it can be seen that the differentiation gives a remainder that is approximately $p$ times bigger. This indicates that the error that will dominate is due to the differentiation of a polynomial that is experiencing aliasing.

%% file: variables.tex
%!TEX root = ./aliasing_main.tex
\section{Primitive and Conserved Variables}
\label{sec:variables}

\subsection{Euler's Equations}
	We will begin by considering the 1D Euler's equations in the conservative form. 
	\begin{equation}\label{eq:euler}
		\px{\mathbf{Q}}{t} + \px{\mathbf{f}(\mathbf{Q})}{x} = 0
	\end{equation}
	for
	\begin{equation}\label{eq:euler_var}
		\mathbf{Q} =\begin{bmatrix}
		\rho \\ \rho u \\ E
		\end{bmatrix}, \quad  \mathrm{and} \quad \mathbf{f}(\mathbf{Q}) =\begin{bmatrix}
		\rho u\\ \rho u^2+p \\ u(E+p)
		\end{bmatrix}.
	\end{equation}
	The concern of this paper is what information should be stored, whilst still solving this equation.

\subsubsection{Conserved Variable Computation}
In an implementation where the conserved variables are not stored directly, if the conservative form of Euler's equations is to be solved, then the conserved variables must, at some stage, be computed.

	\begin{align}\label{eq:p2conv}
		\mathbf{Q}_p &\rightarrow \mathbf{Q}_c \\
		\begin{bmatrix}
			\rho\\u\\p
		\end{bmatrix} &\rightarrow 
		\begin{bmatrix}
			\rho\\\rho u\\ \frac{p}{\gamma-1} + \frac{1}{2}\rho(u^2)
		\end{bmatrix} = 
		\bigo{\begin{bmatrix}\xi^p \\ \xi^{2p} \\ \xi^{3p} \end{bmatrix}}
	\end{align}
	This transformation is shown in Eq.(\ref{eq:p2conv}). It should be clear that if $\mathbf{Q}_p$ is represented by a polynomial of order $p$, then the terms $\rho u$, $\rho v$, and $\rho w$ will be polynomials of order $2p$, while $u(E+p)$ will be of order $3p$. If the order of the scheme is greater than $3p$ this poses no issue. However, depending on how $\mathbf{Q}_c$ is used this could pose a problem.
	
\subsubsection{Inviscid Flux Computation}
	In most implementations seen by the authors, when the primitives are stored they are also subsequently used to form the flux, as opposed to using $\mathbf{Q}_c$. Therefore, the order of the flux variables formed from the primitives is:

	\begin{align}\label{eq:qp2f}
		\mathbf{Q}_p &\rightarrow \mathbf{f} \\
		\begin{bmatrix}
			\rho\\u\\p
		\end{bmatrix}
		&\rightarrow
		\begin{bmatrix}
			\rho u\\\rho u^2 + p\\ u(\frac{\gamma p}{\gamma-1} + \frac{1}	{2}\rho(u^2))
		\end{bmatrix}  = \bigo{\begin{bmatrix}
		\xi^{2p} \\ \xi^{3p} \\ \xi^{4p} \end{bmatrix}}
	\end{align}

If instead the conserved variables are used we obtain:
	\begin{align}\label{eq:qc2f}
		\mathbf{Q}_c &\rightarrow \mathbf{f} \\
		\begin{bmatrix}
			\rho\\ \rho u\\E
	\end{bmatrix}
	&\rightarrow
	\begin{bmatrix}
		(\rho u)\\ \frac{(\rho u)^2}{\rho} + (\gamma-1)\Big(E - \half\frac{(\rho u)^2}{\rho}\Big)\\ \frac{(\rho u)}{\rho}\Big(\gamma E - \half(\gamma-1)\frac{(\rho u)^2}{\rho}\Big)
	\end{bmatrix} = 
	\bigo{\begin{bmatrix} \xi^{p} \\ \xi^{2p}/\xi^p \\ \xi^{3p}/\xi^{2p} \end{bmatrix}}
	\end{align}

To clarify the notation used here, it is intended for $\bigob{\xi^{2p}/\xi^{p}}$ to mean a $2p^{\mathrm{th}}$ order polynomial divided by a $p^{\mathrm{th}}$ order polynomial. If the polynomial $1/\bigob{\xi^p}$ is then expanded about zero to form a series of monomials, the series is $\bigob{\xi^\infty}$. From this, it can be seen that in Eq.(\ref{eq:qc2f}), we have avoided the $\xi^{4p}$ term, but at the expense of dividing by $\rho$. This raises the question as to whether this formulation is more accurate --- specifically, is the convergence of the $1/\rho$ series sufficiently fast to reduce aliasing? Importantly though, this method avoids a whole mechanism of the aliasing, introduced through the conversion in Eq.(\ref{eq:p2conv}).

	Another option that will be explored is storing the conserved variables, but with energy substituted for pressure, $\mathbf{Q}_{c+p}$. The reason being that in industrial codes pressure is used frequently and this option would reduce the work involved in converting an implementation. Hence, the conversion from $\mathbf{Q}_{c+p}$ to the flux, $\mathbf{f}$, is:
	\begin{align}\label{eq:qcp2f}
		\mathbf{Q}_{c+p} &\rightarrow \mathbf{f} \\
		\begin{bmatrix}
			\rho\\\rho u \\p
		\end{bmatrix}
		&\rightarrow
		\begin{bmatrix}
			(\rho u)\\ \frac{(\rho u)^2}{\rho} + p\\ \frac{(\rho u)}{\rho}\Big(\frac{\gamma p }{\gamma -1} + \half\frac{(\rho u)^2}{\rho}\Big)
	\end{bmatrix} = 
	\bigo{\begin{bmatrix} \xi^{p} \\ \xi^{2p}/\xi^p \\ \xi^{3p}/\xi^{2p} \end{bmatrix}}
	\end{align}
	This method will also require a conversion step to retrieve the conserved variables if Eqs.(\ref{eq:euler}~\&~\ref{eq:euler_var}) are to be solved. This then introduces aliasing of order:
	\begin{align}
		\mathbf{Q}_{c+p} &\rightarrow \mathbf{Q}_c \\
		\begin{bmatrix}
		\rho \\\rho u \\p
		\end{bmatrix} &\rightarrow
		\begin{bmatrix}
			(\rho) \\ (\rho u) \\ \frac{p}{\gamma-1} + \half\frac{(\rho u)^2}{\rho}
		\end{bmatrix} = 
		\bigo{\begin{bmatrix} \xi^p \\ \xi^p \\ \xi^{2p}/\xi^p \end{bmatrix}}
	\end{align}	
	This method has the potential to reduce the aliasing in forming the conserved variables and flux, as there is no longer the $\xi^{4p}$ that is present in Eq.(\ref{eq:qp2f}). However, this is again dependent on the nature of $1/\rho$. 

\subsection{Navier-Stokes Equations}
	To confront more complex problems of fluid dynamical relevance, it is essential to consider the Navier-Stokes equations, written for 3D in the conservative form as:

	\begin{equation}
		\px{\mathbf{Q}}{t} + \nabla\cdot\mathbf{F}(\mathbf{Q},\nabla\mathbf{Q}) = 0
	\end{equation}
	where 
	\begin{equation}
		\nabla\cdot\mathbf{F} = 
			(\mathbf{f}^\mathrm{inv}-\mathbf{f}^\mathrm{vis})_x +
			(\mathbf{g}^\mathrm{inv}-\mathbf{g}^\mathrm{vis})_y +
			(\mathbf{h}^\mathrm{inv}-\mathbf{h}^\mathrm{vis})_z
	\end{equation}

	If we take the bulk viscosity, $\mu_b$, to be zero, then $\mathbf{f}^\mathrm{vis}$ can be defined as:
\begin{equation}
	\mu\begin{bmatrix} 
		0 \\ \tau_{xx} \\ \tau_{xy} \\ \tau_{xz} \\ u\tau_{xx} + v\tau_{xy} + w\tau_{xz} + \frac{\kappa}{\mu}T_x
	\end{bmatrix} = \mu\begin{bmatrix} 
		0 \\ \frac{4}{3}u_x - \frac{2}{3}(v_y + w_z) \\ u_y + v_x \\ w_x + u_z \\ u(\frac{4}{3}u_x - \frac{2}{3}(v_y + w_z)) + v(u_y + v_x) + w(w_x + u_z) + \frac{\kappa}{\mu}T_x
	\end{bmatrix}
\end{equation}
	with $\mathbf{g}^\mathrm{vis}$ and $\mathbf{h}^\mathrm{vis}$ similarly defined.
	
	The importance of considering this equation is that, due to phenomena such as the energy cascade, in a method which does not suffer from implementation aliasing, aliasing will arise in LES due to the partial resolution of vortical motions. Hence, for turbulent flows, any difference is likely to be more marked as implementation aliasing amplifies the existing numerical aliasing.

    Clearly for the case when primitive variables are stored, the gradients of the primitive can be directly calculated and used to form the viscous flux. However, when the conserved variables are stored there are two options available to form the gradients needed: to convert the conserved variables to the primitives and to calculate the gradients needed directly:
	\begin{equation}
		\begin{bmatrix}
			\rho\\ \rho u\\ \rho v\\ \rho w\\E
		\end{bmatrix} \rightarrow
		\begin{bmatrix}
			\rho \\ u \\ v \\ w \\ p
		\end{bmatrix} \rightarrow
		\begin{bmatrix}
		\rho_x & \dots\\ 
		u_x & \dots\\
		v_x & \dots\\
		w_x & \dots\\
		\frac{c_v}{\gamma-1}(\rho^{-1}p_x-\rho^{-2}p\rho_x) & \dots
		\end{bmatrix};
	\end{equation}
	or to calculate the gradient of the conserved variables and use the product rule to convert them to what is needed:
	\begin{equation}
		\begin{bmatrix}
		\rho_x & \rho_y & \rho_z \\
		(\rho u)_x & (\rho u)_y & (\rho u)_z \\
		(\rho v)_x & (\rho v)_y & (\rho v)_z \\
		(\rho w)_x & (\rho w)_y & (\rho w)_z \\
		E_x & E_y & E_z
		\end{bmatrix} \rightarrow
		\begin{bmatrix}
		\rho_x & \rho_y & \rho_z \\
		u_x & u_y & u_z \\
		v_x & v_y & v_z \\
		w_x & w_y & w_z \\
		T_x & T_y & T_z
		\end{bmatrix}.
	\end{equation}
	These two options can be simplified as:
	\begin{alignat}{6}
		&\mathbf{Q}_c &&\rightarrow &&\mathbf{Q}_p  &&\rightarrow &&\nabla\mathbf{Q}_p \label{eq:cons} \\
		&\mathbf{Q}_c &&\rightarrow \nabla&&\mathbf{Q}_c &&\rightarrow &&\nabla\mathbf{Q}_p \label{eq:cong}
	\end{alignat}
	where $\nabla\mathbf{Q}$ is the gradient of $\mathbf{Q}$. Here the final row of $\nabla\mathbf{Q}$ is the gradient of temperature, $\nabla T$, for convenience in the calculation of the viscous flux.
	
	The method for calculating the required gradients from the product rule applied to the conserved variable gradient formulation is:
	\begin{equation}
		\frac{1}{\rho}\begin{bmatrix}
		\rho\rho_x & \dots\\
		\big((\rho u)_x - \rho^{-1}(\rho u)\rho_x\big) & \dots\\
		\big((\rho v)_x - \rho^{-1}(\rho v)\rho_x\big) & \dots\\
		\big((\rho w)_x - \rho^{-1}(\rho w)\rho_x\big) & \dots\\
		\big(E_x - \rho^{-1}E\rho_x\big) - \big((\rho u)u_x + (\rho v)v_x + (\rho w)w_x  \big) & \dots
		\end{bmatrix} = \begin{bmatrix}
		\rho_x & \rho_y & \rho_z \\
		u_x & u_y & u_z \\
		v_x & v_y & v_z \\
		w_x & w_y & w_z \\
		T_x & T_y & T_z
		\end{bmatrix}
	\end{equation}		
	The polynomial order of this step is then: 
	\begin{equation}
		\begin{split}
		\frac{1}{\rho}\begin{bmatrix}
		\rho\rho_x\\
		\big((\rho u)_x - \rho^{-1}(\rho u)\rho_x\big)\\
		\big((\rho v)_x - \rho^{-1}(\rho v)\rho_x\big)\\
		\big((\rho w)_x - \rho^{-1}(\rho w)\rho_x\big)\\
		\big(E_x - \rho^{-1}E\rho_x\big) - \big((\rho u)u_x + (\rho v)v_x + (\rho w)w_x  \big)
		\end{bmatrix} = & \\
		\bigo \begin{bmatrix}
		\xi^{p-1}(\eta\zeta)^p \\
		\xi^{p-1}(\eta\zeta)^p/(\xi\eta\zeta)^{p} + \xi^{2p-1}(\eta\zeta)^{2p}/(\xi\eta\zeta)^{2p} \\
		\xi^{p-1}(\eta\zeta)^p/(\xi\eta\zeta)^{p} + \xi^{2p-1}(\eta\zeta)^{2p}/(\xi\eta\zeta)^{2p} \\
		\xi^{p-1}(\eta\zeta)^p/(\xi\eta\zeta)^{p} + \xi^{2p-1}(\eta\zeta)^{2p}/(\xi\eta\zeta)^{2p} \\
		\xi^{p-1}(\eta\zeta)^p/(\xi\eta\zeta)^{p} + \xi^{2p-1}(\eta\zeta)^{2p}/(\xi\eta\zeta)^{2p} + \xi^{2p-1}(\eta\zeta)^{2p}/(\xi\eta\zeta)^p\\
		\end{bmatrix} &
		\end{split}
	\end{equation}
	Again, it should be clear that the momentum and energy (rows 2-5) terms experience the most aliasing, although it is not clear what effect that the division will have on aliasing. However, it is likely that the decay rate of the infinite quotient series will be fast in most cases. 
	
	The methods of data storage that will be investigated for the Navier-Stokes equations and Euler's equations (where applicable) can then be summerised as storing the:

	\begin{itemize}
		\item[$\bullet$ (\emph{A})] Primitive variables
		\item[$\bullet$ \quad \;\;] Conserved variables, with the gradients for the Navier-Stokes equation calculated from the:  
		\begin{itemize}
			\item[-- (\emph{B})] conserved variables converted to primitive variables
			\item[-- (\emph{C})] product rule applied to the gradient of the conserved variables 
		\end{itemize}
		\item[$\bullet$ (\emph{D})] Conserved variables, but with pressure instead of energy
	\end{itemize}
	where the letters in brackets are shorthand identifiers that will be used in the subsequent commentary.
	
	At this point, we wish to link the ideas presented in section~\ref{sec:palias} with the methods of this section. It should be clear that in order for this form of aliasing error to be incorporated into the solution, then at some stage interpolation or polynomial fitting has to be used within the calculation. For FR, this comes when the gradient is calculated or the edge points are extrapolated from the points inside the element. However, if only the nodal values are used, as is the case in second-order Finite Volume (FV) methods, then there is no mechanism by which this form of aliasing can affect the solution. Take the example of converting primitive variables to conservative variables, and back again:
	\begin{equation}
		\mathbf{Q}_p \rightarrow \mathbf{Q}_c \rightarrow \mathbf{Q}^\prime_p.
	\end{equation}
	It should be apparent that beyond rounding error introduced, $\mathbf{Q}_p = \mathbf{Q}^\prime_p$. Therefore, the means of variable storage will not affect FV but will affect any method that in some way interpolates or fits a polynomial.

%% file: icv.tex
%!TEX root = ./aliasing_main.tex
\section{Isentropic Convecting Vortex}
\label{sec:icv}
	To evaluate the impact of the changes suggested in section~\ref{sec:variables} we will begin by studying the effect on the error and total kinetic energy of the Isentropic Convecting Vortex (ICV)~\cite{Shu1997}. The ICV is of interest as it is an analytical solution to Euler's equations and hence allows for the error at a given time to be calculated. The problem that we are confronted with when using high order, the ICV, and a periodic domain, is that the solution is only guaranteed to be $C^0$ continuous. This can be understood by considering the initial condition:
	\begin{align}
		\rho &= \bigg(1-\frac{(\gamma-1)\beta^2}{8\gamma\pi^2}\exp{(1-r^2)}\bigg)^{\frac{1}{\gamma-1}} \\
		u &= u_0 + \frac{\beta}{2\pi}(y_0-y)\exp{\bigg(\frac{1-r^2}{2}\bigg)}\\
		v &= v_0 + \frac{\beta}{2\pi}(x-x_0)\exp{\bigg(\frac{1-r^2}{2}\bigg)}\\
		w &= 0 \\
		p &= \bigg(1-\frac{(\gamma-1)\beta^2}{8\gamma\pi^2}\exp{(1-r^2)}\bigg)^{\frac{\gamma}{\gamma-1}} \\
		r^2 & = (x-x_0)^2 + (y-y_0)^2
	\end{align}  
	where $u_0$ and $v_0$ are the advective velocities and $\beta$ is the vortex strength (typically $\beta=5$ is used). Hence, it can be seen that as the distance $r$ is increased the vortex slowly decays and, on a finite but periodic domain, this will lead to discontinuities in the gradient. This is a point that will be of importance later when reviewing results.
	
	The metrics that we will use to review the accuracy are the point averaged absolute error in the density:
	\begin{equation}
		e(t) = \frac{1}{N_p}\sum^{N_p}_{i=1}|\rho_i-\rho(\mathbf{x}_i,t)|_2	
	\end{equation}  
	and the total kinetic energy:
	\begin{equation}
		E_k(t) = \frac{1}{2|\mathbf{\Omega}|} \int_{\mathbf{\Omega}}\rho\mathbf{V}\cdot\mathbf{V}\mathrm{d}\mathbf{x}
	\end{equation}
	where $|\mathbf{\Omega}|$ is the domain volume.	
	
	\begin{figure}
		\centering
			\includegraphics[width=0.5\linewidth]{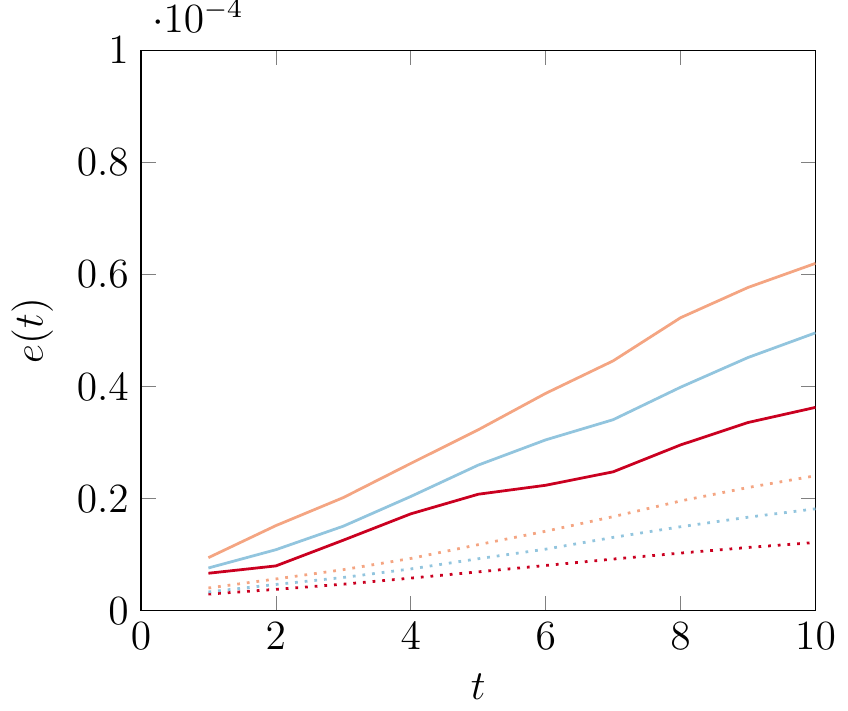}
			\\
			\includegraphics[width=0.55\linewidth]{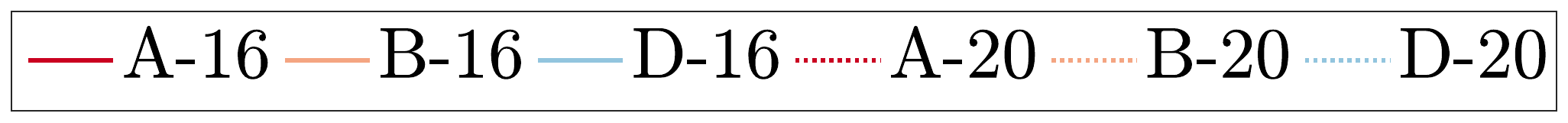}
		\caption{Variation of error in ICV density with time for FR, $p=4$, using methods A, B and D on $16\times16\times2$ and $20\times20\times2$ element grids.}
		\label{fig:FRp4_ICVe}
	\end{figure}
	
	We begin by investigating the effect of storing the primitive variables (A), the conserved variables (B), and the conserved variables with energy substituted for pressure (D), on the error. This is shown in Fig.~\ref{fig:FRp4_ICVe}. Clearly, method A has the lowest levels of error followed by D then B, and this ordering does not change as the grid is refined. This results may be thought to be contrary to the expected outcome. However, to appreciate what is going on, consider the development of the kinetic energy with time.
	
	\begin{figure}
		\centering		
		\begin{subfigure}[b]{0.48\linewidth}
			\centering
			\includegraphics[width=\linewidth]{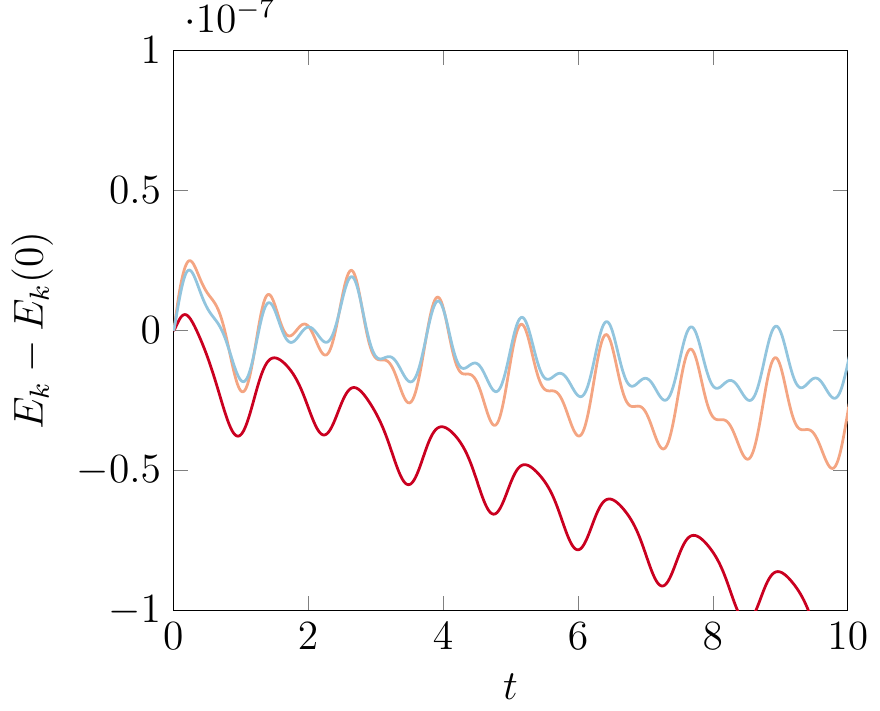}
			\caption{$16\times16\times2$ elements}
			\label{fig:FRp4_ICV_energy_16}
		\end{subfigure}
		~
		\begin{subfigure}[b]{0.48\linewidth}
			\centering
			\includegraphics[width=\linewidth]{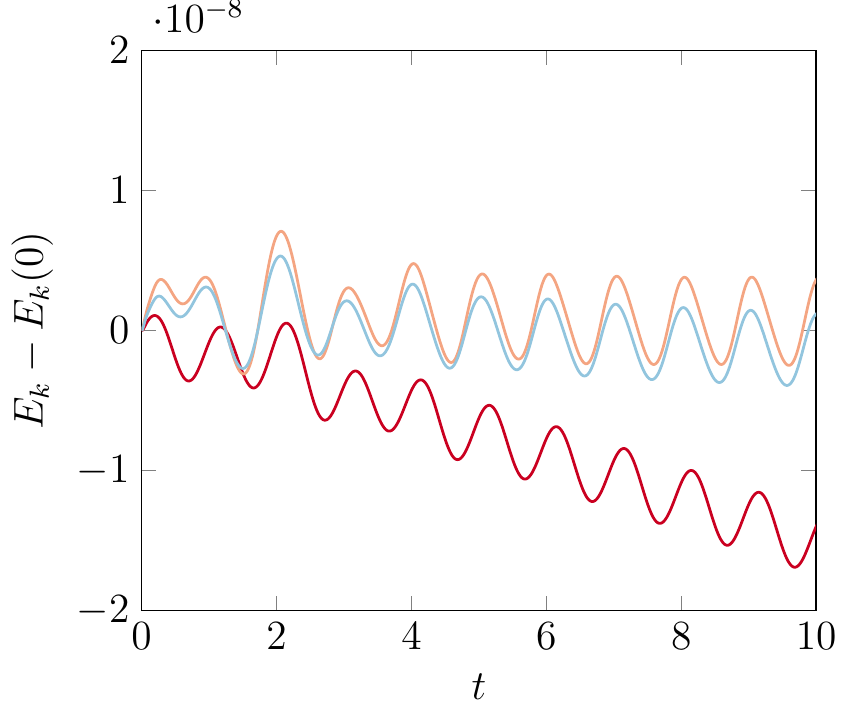}
			\caption{$20\times20\times2$ elements}
			\label{fig:FRp4_ICV_energy_20}
		\end{subfigure}
		~
		\begin{subfigure}[b]{0.18\linewidth}
			\centering
			\includegraphics[width=\linewidth]{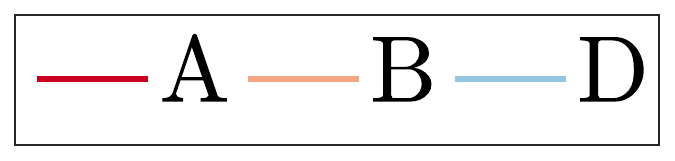}
		\end{subfigure}
		\caption{Variation in total kinetic energy of the ICV, FR $p=4$, for two grid resolutions. using methods A, B, and D.}
		\label{fig:FR_ICV_energy}
	\end{figure}
	
	Figure~\ref{fig:FR_ICV_energy} shows how the kinetic energy in the domain changes with time. For both grid resolutions, the rate of kinetic energy dissipation of A is higher than B and D, while B and D are similar. In the higher resolution case, Fig.~\ref{fig:FRp4_ICV_energy_20}, the dissipation of B is found to be less than that of D. The importance of this is that, as was stated earlier, the ICV initial condition is only formally $C^0$ continuous. Therefore, the lower dissipation that methods B and D exhibit leads to the errors introduced via the discontinuities in the gradient not being as damped as in the case of A. Hence, the error grows faster while also showing less dissipation. Ultimately, it looks as though methods B and D can aid in the reduction of dissipation via aliasing, although this does have some associated issues.

%% file: tgv.tex
%!TEX root = ./aliasing_main.tex
\section{Taylor-Green Vortex}
\label{sec:tgv}

	The final investigation to be considered is the application of the various forms of stored variable to the full Navier-Stokes equations for a flow which exhibits turbulence. The flow of choice for this is the canonical Taylor-Green vortex~\cite{Taylor1937}, where the exact flow field used is defined by DeBonis~\cite{DeBonis2013,Chapelier2012}. This case is chosen as not only is it a case for the Navier-Stoke equations, but it exhibits a transition from an inviscid regime to a fully turbulent flow, via the mechanism of vortex stretching and shearing. This is key as, not only is it more representative of real engineering flows, but a transition to turbulence will introduce an energy cascade to the flow and hence induce aliasing.

	The key non-dimensional parameters used in the characterisation of the physics of this flow are:
	\begin{equation}
	R_e = \frac{\rho_0U_0L}{\mu}, \quad P_r = \frac{\mu\gamma R}{\kappa(\gamma-1)}, \quad M_a = \frac{U_0}{\sqrt{\gamma RT_0}}.
	\end{equation}

	Here the free parameters that define the initial condition of the flowfield are: the stagnation density, $\rho_0$, the stagnation pressure, $p_0$, the stagnation temperature, $T_0$, and the velocity magnitude, $U_0$. For all the TGV tests, the Prandtl number was set at $P_r=0.71$, with a bulk viscosity of zero. The metrics that we will use to study the behaviour of the numerical method applied to the TGV are the rate of kinetic energy dissipation and enstrophy dissipation:
	\begin{align}
		\epsilon_1 &= -\dx{E_k}{t} = -\dx{}{t}\bigg(\frac{1}{2\rho_0U_0^2|\mathbf{\Omega}|} \int_{\mathbf{\Omega}}\rho\mathbf{V}\cdot\mathbf{V}\mathrm{d}\mathbf{x}\bigg)\\
		\epsilon_2 &= \frac{\mu}{\rho_0^2U_0^2|\mathbf{\Omega}|}\int_{\mathbf{\Omega}}\rho(\pmb{\omega}\cdot\pmb{\omega})\mathrm{d}\mathbf{x}
	\end{align}
	where $\pmb{\omega}$ is the vector of vorticity, $\mu$ is the shear viscosity, and where $\epsilon_1$~\&~$\epsilon_2$ have been normalised. 
	
	 Aliasing is the main focus of this paper, and as such, we want to investigate if the different method of variable storage impacts the accuracy of the solution. As a result, there are two things which will be varied, the first of which the Reynolds number. Here, Reynolds numbers of $R_e=400,1600,3000$, are used, with reference DNS data (denoted by ``ref'') available from \cite{Brachet1983}. This is because it will trigger a variety of different physics. The second variable that is varied is the Mach number, where values of $M_a=0.08$, and $0.31$ are used. The effect of compressibility on the TGV was investigated by~\cite{Peng2018} at various Mach numbers between $0.5$ and $2$, with $0.5$ not being found to exhibit shocklets. Therefore, testing at $M_a=0.31$ will test the introduction of aliasing due to larger spatial variations in $\rho$, but without triggering issues relating to shock capturing.
	
	For the majority of the investigation, a 3D Navier-Stokes FR scheme will be used. The grid topology used will be hexahedral, constructed using a tensor product construction of the 1D FR scheme. More details on this construction of FR can be found in \cite{Witherden2014,Castonguay2011,Williams2014} including the method of extension to diffusion equations. The method of calculating the inviscid common interface flux chosen is a Rusanov flux with Davis wave speeds \cite{Rusanov1961,Davis1988}. The viscous common interface flux is found using Bassi~and~Rebays' BR1 scheme \cite{Bassi1997a,Bassi1998a}. The aim of this paper is not to explore the effect of FR correction function, but aliasing. Because of this, an FR correction function that recovers Nodal DG is used throughout~\cite{Huynh2007,Hesthaven2008}.
	
	\begin{figure}
		\centering
		\begin{subfigure}[b]{0.48\linewidth}
			\centering
			\includegraphics[width=\linewidth]{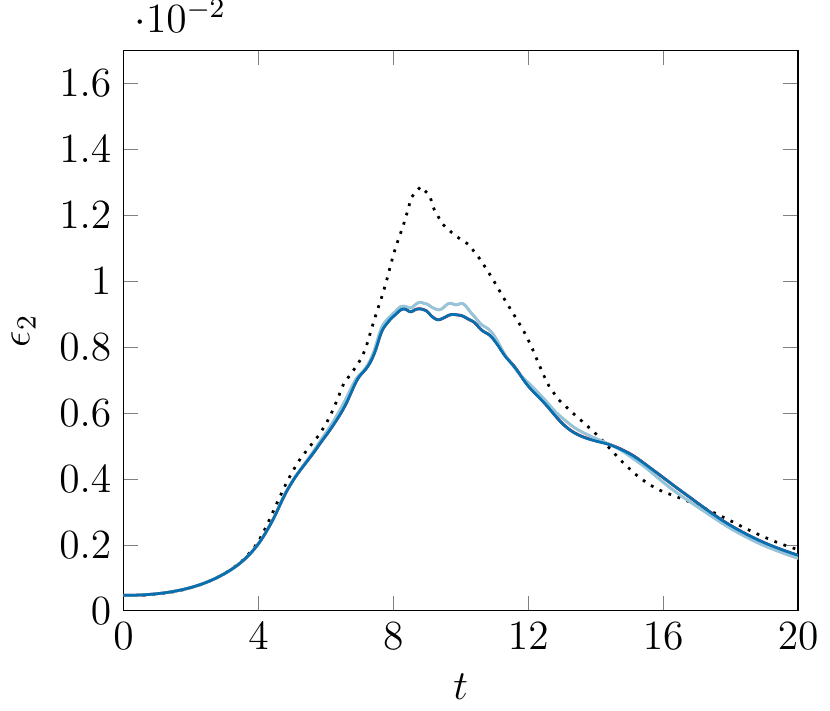}
			\caption{$M_a=0.08$}
			\label{fig:FRp4_1600_m08}
		\end{subfigure}
		~
		\begin{subfigure}[b]{0.48\linewidth}
			\centering
			\includegraphics[width=\linewidth]{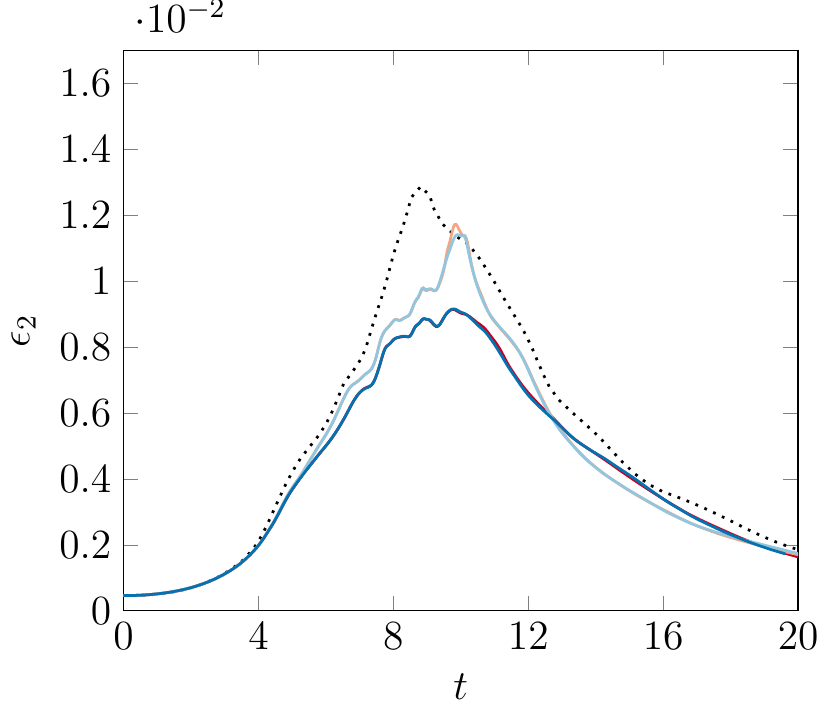}
			\caption{$M_a=0.31$}
			\label{fig:FRp4_1600_m31}
		\end{subfigure}
		~
		\begin{subfigure}[b]{0.33\linewidth}
			\centering
			\includegraphics[width=\linewidth]{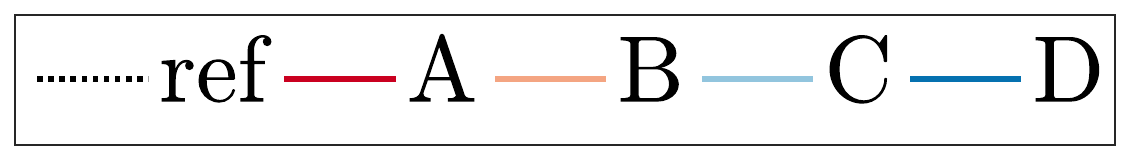}
		\end{subfigure}
		\caption{Enstrophy of the Taylor-Green Vortex with $R_e=1600$, $p=4$ and $80^3$ degrees of freedom for storage methods A-D.}
		\label{fig:FR_1600}
	\end{figure}	
	
	Let us first consider the TGV case when $R_e=1600$ at Mach numbers $0.08$ and $0.31$. We will look to compare the four methods presented in section~\ref{sec:variables} for FR with $p=4$ on a mesh with $80^3$ degrees of freedom, the results of which are shown in Fig.~\ref{fig:FR_1600}. It is seen that when $M_a=0.08$, Fig.~\ref{fig:FRp4_1600_m08}, there is a small increase in the enstrophy disipation when storing data as conserved variables over primitive variables. The results of the conservative variables with pressure instead of energy can be seen to be almost identical to the primitive variable results. It can also be noted that the largest difference is seen around the time of peak dissipation, and not in the region $4<t<7$. This seems to indicate that the effect of changing the method of variable storage is to reduce the numerical/aliasing based dissipation at the smallest scales. It is apparent that it does not introduce extra sources of dispersion which would cause excess dissipation around $4<t<7$, when small scales begin to develop in the flow.
	
	Moving on to the case when $M_a=0.31$, the high Mach number will introduce larger spatial variation in the density as the flow becomes more compressible. This should increase the effect of aliasing error on the solution. The enstrophy is displayed in Fig.~\ref{fig:FRp4_1600_m31} and clearly shows a far larger change between the full conservative and the primitive methods. Again, the cases of primitive and partial conservative with pressure are similar, this indicates that the improvement is largely originating from the change in the handling of the energy equation. The work in section~\ref{sec:variables} showed that when the primitive data is stored, the formulation of energy from pressure see aliasing error $p$ orders higher than the momentum terms, and hence it is to be expected that the largest contribution to the improving the scheme comes from the energy equation. 
	
	At the higher Mach number, there is a noticeable difference in Fig.~\ref{fig:FRp4_1600_m31} between the fully conservative with the gradient calculated from the converted primitives and the gradient calculated from the application of the product rule. It is hard to attribute this difference to a particular aspect, but we will explore this further. 
	
	\begin{figure}
		\centering
		\begin{subfigure}[b]{0.48\linewidth}
			\centering
			\includegraphics[width=\linewidth]{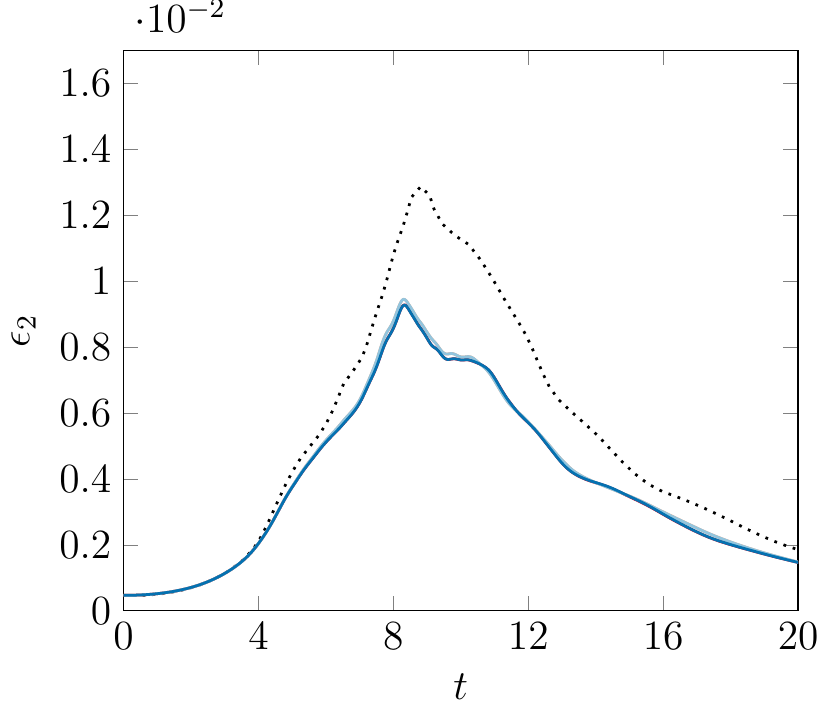}
			\caption{$M_a=0.08$}
			\label{fig:FRp3_1600_m08}
		\end{subfigure}
		~
		\begin{subfigure}[b]{0.48\linewidth}
			\centering
			\includegraphics[width=\linewidth]{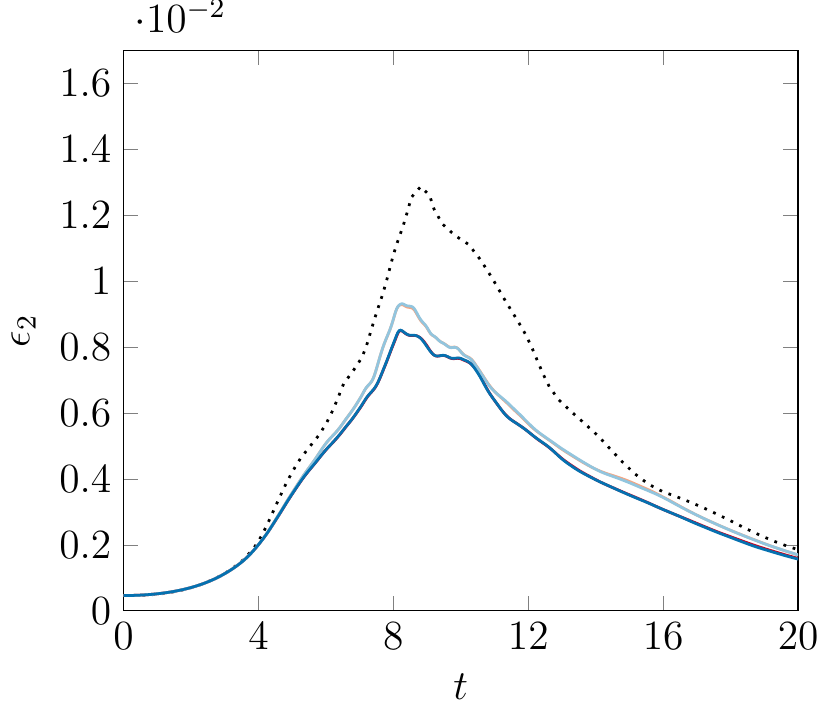}
			\caption{$M_a=0.31$}
			\label{fig:FRp3_1600_m31}
		\end{subfigure}
		~
		\begin{subfigure}[b]{0.33\linewidth}
			\centering
			\includegraphics[width=\linewidth]{TGV_leg.pdf}
		\end{subfigure}
		\caption{Enstrophy of the Taylor-Green Vortex with $R_e=1600$, $p=3$ and $80^3$ degrees of freedom.}
		\label{fig:FR3_1600}
	\end{figure}
	
	In section~\ref{sec:palias} the dependency of interpolation rounding error on order, and its factorial increase with order was shown analytically. To investigate the effect of order we consider the case of $R_e=1600$ run at $p=3$ for the same number of degrees of freedom. The results of this are shown in Fig.~\ref{fig:FR3_1600}. By comparison of Figs.~\ref{fig:FRp3_1600_m08}~\&~\ref{fig:FRp3_1600_m31}, it can be seen that there is still a larger difference between the methods in the high Mach number case than at low Mach number. However, when comparing Figs.~\ref{fig:FR3_1600}~\&~\ref{fig:FR_1600}, the difference between methods is markedly smaller at lower order. This evidence is in agreement with the earlier analytical predictions. As we move to a higher order, this mechanism of aliasing becomes increasingly important.
	
	\begin{figure}
		\centering
		\begin{subfigure}[b]{0.48\linewidth}
			\centering
			\includegraphics[width=\linewidth]{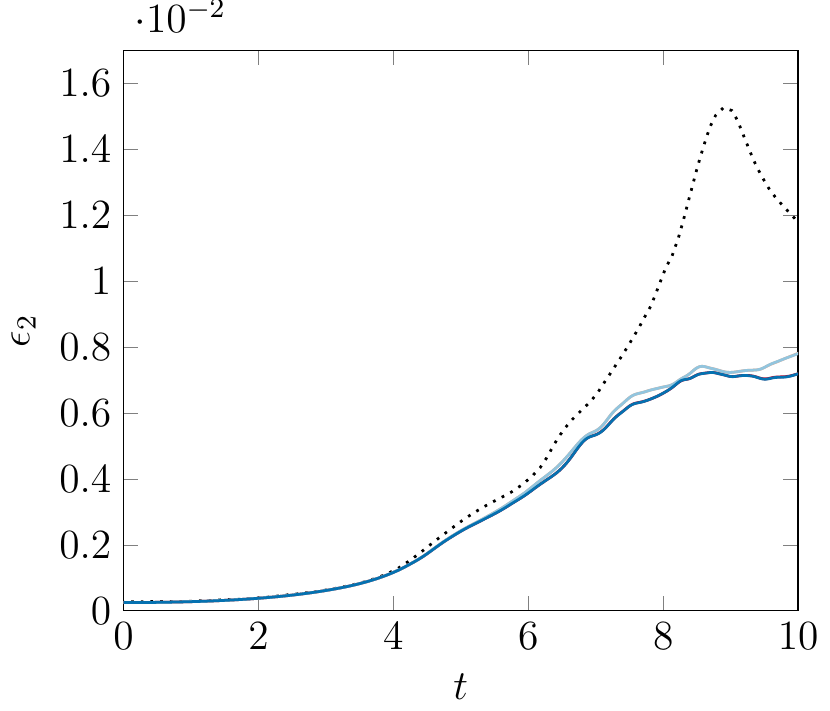}
			\caption{$M_a=0.08$}
			\label{fig:FRp4_3000_m08}
		\end{subfigure}
		~
		\begin{subfigure}[b]{0.48\linewidth}
			\centering
			\includegraphics[width=\linewidth]{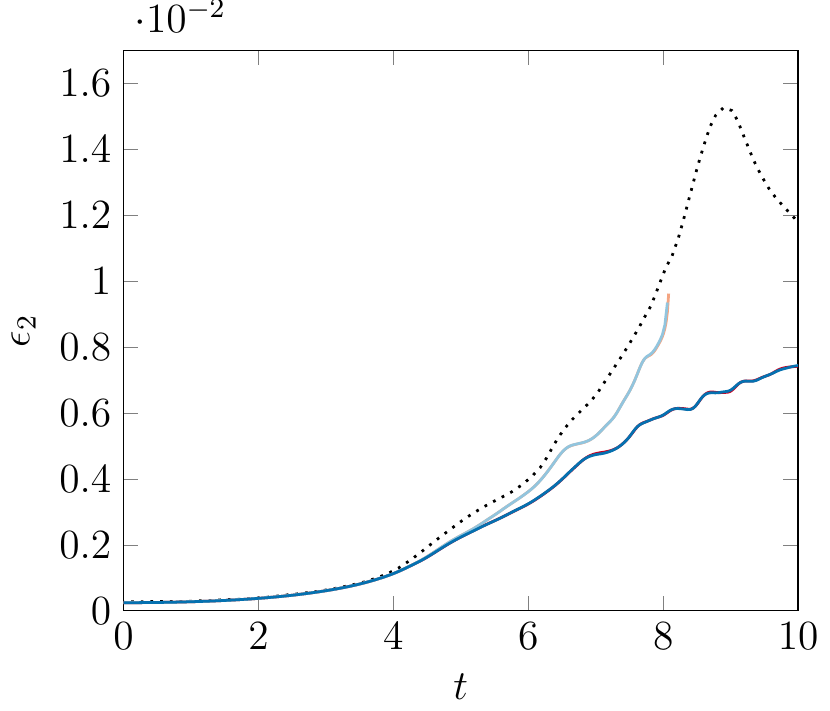}
			\caption{$M_a=0.31$}
			\label{fig:FRp4_3000_m31}
		\end{subfigure}
		~
		\begin{subfigure}[b]{0.33\linewidth}
			\centering
			\includegraphics[width=\linewidth]{TGV_leg.pdf}
		\end{subfigure}
		\caption{Enstrophy of the Taylor-Green Vortex with $R_e=3000$, $p=4$ and $80^3$ degrees of freedom.}
		\label{fig:FR_3000}
	\end{figure}
	
	We will now explore the effect of increasing the Reynolds numbers for the same grid resolution. In particular, we choose $R_e=3000$, which was explored with DNS by Brachet~\etal~\cite{Brachet1983} and with DG by Chapelier~\etal~\cite{Chapelier2012}. The results of the application of $p=4$ FR with the various methods of storage are presented in Fig.~\ref{fig:FR_3000}. Firstly studying the $M_a=0.08$ case, there is again a noticeable difference between the conservative and primitive enstrophy, whichcan be attributed to the decrease in numerical/aliasing based dissipation, due to the absence of over dissipation when small scales begin to be generated and the increase in dissipation at the expected peak. Hence, the small scales are being preserved for longer thus enabling their increased contribution to physical dissipation.
	
	When the Mach number is increased to $M_a=0.31$ we initially see a larger difference between the formulations, followed by the solution diverging. A similar divergence was observed by Chapelier~\etal~\cite{Chapelier2012} when using DG on an under-resolved mesh. They attributed the divergence to insufficient numerical dissipation to stabilise the under-resolved grid. This adds weight to the argument that the change in the variables stored is mainly reducing numerical dissipation by reducing aliasing and not introducing dispersion. It should be noted that our results differ slightly from those of Chapelier~\etal~\cite{Chapelier2012} as they are solving the filtered LES equations, whereas we are using implicit LES.
	
	\begin{figure}
		\centering
		\begin{subfigure}[b]{0.48\linewidth}
			\centering
			\includegraphics[width=\linewidth]{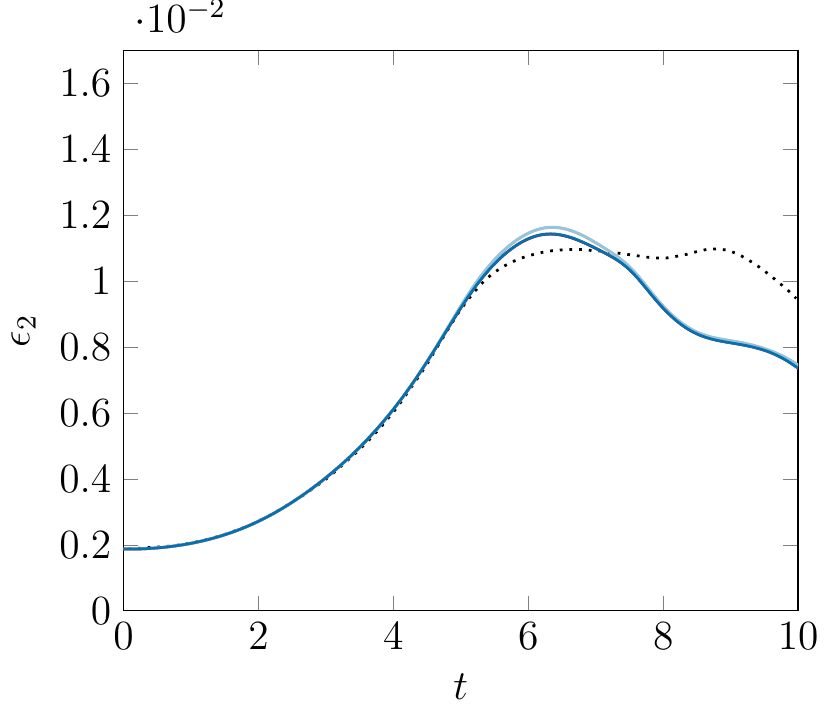}
			\caption{$M_a=0.08$}
			\label{fig:FRp4_400_m08}
		\end{subfigure}
		~
		\begin{subfigure}[b]{0.48\linewidth}
			\centering
			\includegraphics[width=\linewidth]{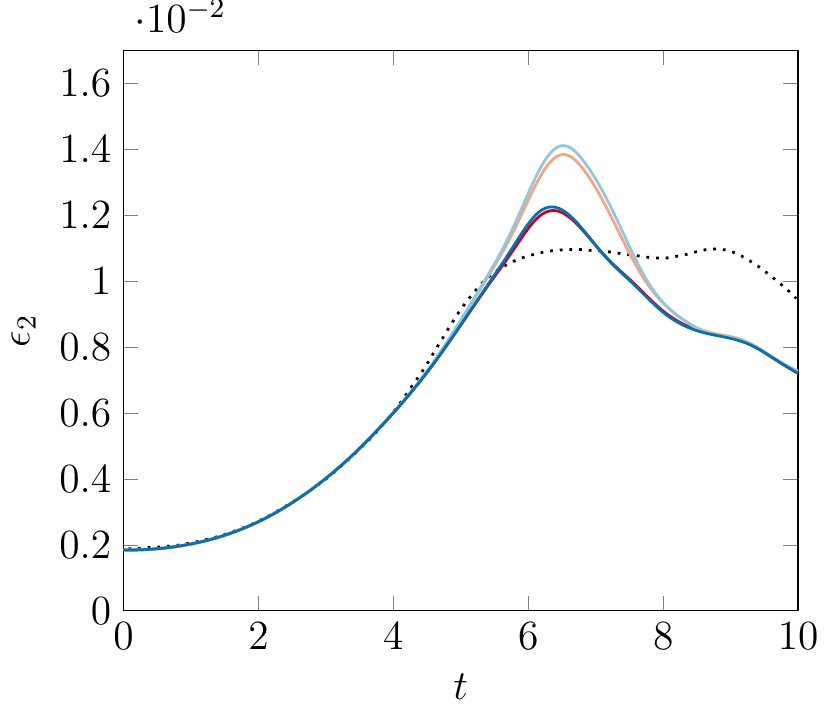}
			\caption{$M_a=0.31$}
			\label{fig:FRp4_400_m31}
		\end{subfigure}
		~
		\begin{subfigure}[b]{0.33\linewidth}
			\centering
			\includegraphics[width=\linewidth]{TGV_leg.pdf}
		\end{subfigure}
		\caption{Enstrophy of the Taylor-Green Vortex with $R_e=400$, $p=4$ and $40^3$ degrees of freedom.}
		\label{fig:FR_400}
	\end{figure}
	
	To highlight the impact of using the conserved variables with the product rule to calculate the gradient of the primitives, Eq.(\ref{eq:cong}), we will reduce the Reynolds number to $R_e=400$, such that the viscous terms become more important. Through testing, it was found that when $80^3$ degrees of freedom were used the case was highly resolved, with the enstrophy based decay rate lying on top of the DNS results. Therefore, to introduce a source of aliasing, the grid resolution was reduced to a level that is more in keeping with LES --- in this case to a $R_{e,\mathrm{cell}} = 50$ at $p=4$, or $40^3$ degrees of freedom, the results of which are shown in Fig.~\ref{fig:FR_400}.
	
%	\begin{figure}
%		\centering
%		\captionof{table}{}
%		\begin{tabular}{|c|c|}	
%			\hline
%			Type & Runtime (ns) \\ \hline
%			A & 6422690 \\ \hline
%			B & 5831547 \\ \hline
%			C & 5634136 \\ \hline
%			D & 6756803 \\ \hline
%		\end{tabular}
%		\label{tab:runtime_comp}
%	\end{figure}
	\begin{figure}
		\centering
		\captionof{table}{Computation time comparison for one full RK44 explicit time step on a $8^3$, p=4, mesh. Time saving shown relative to scheme A.}
		\begin{tabular}{|c|c|c|}
			\hline
			Type & Computation time (ms) & Time Saving ($\%$)\\ \hline
			A & 6.423 & -- \\ \hline
			B & 5.832 & $9.2$ \\ \hline
			C & 5.634 & $12.3$ \\ \hline
			D & 6.757 & $-5.2$ \\ \hline
		\end{tabular}
		\label{tab:runtime_comp}
	\end{figure}
	
	The different methods outlined in section~\ref{sec:variables} will obviously require differing numbers of floating point operations as some conversion steps are required or different numbers of multiplications to build things such as the flux terms. Therefore, we wish to understand what the impact on computational performance is, to this end we will profile the implementation. The implementation of FR used is an in house FR solver called Forflux, written in Fortran with Cuda Fortran and cuBLAS, both version 9.1, for GPU acceleration. The current implementation for small cases leads to the entire memory space being resident on the GPU and hence the CPU plays little to no role in the computation. The case profiled is a TGV, $p=4$, with $8^3$ elements run on a Titan Xp. Using the profiler, \texttt{pgprof}, the runtime for one complete explicit time step was found and is shown in Table~\ref{tab:runtime_comp}.
	
	It is clear that the continual conversion to or from the primitive variables has a noticeable impact on the computational time. In this case, the method that required the fewest number of conversions, method C (conservative variables using the product rule to calculate the gradient of the primitives), was the fastest. Method C gave a $12.3\%$ reduction in computational time, which, all other things being equal, makes this a reasonable optimisation strategy to consider. Method D on the other hand, (conservative variables with $E$ swapped for $p$), was slower as there are even more conversions required than the baseline primitive method.

	\begin{figure}
		\centering
		\begin{subfigure}[b]{0.48\linewidth}
			\centering
			\includegraphics[width=\linewidth]{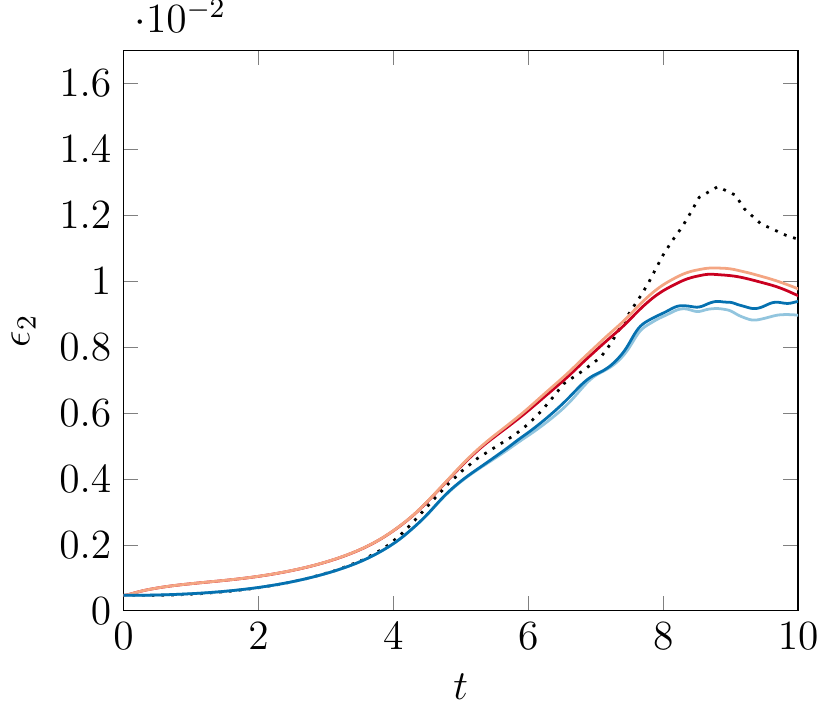}
			\caption{$M_a=0.08$}
			\label{fig:FRp4_fp_1600_m08}
		\end{subfigure}
		~
		\begin{subfigure}[b]{0.48\linewidth}
			\centering
			\includegraphics[width=\linewidth]{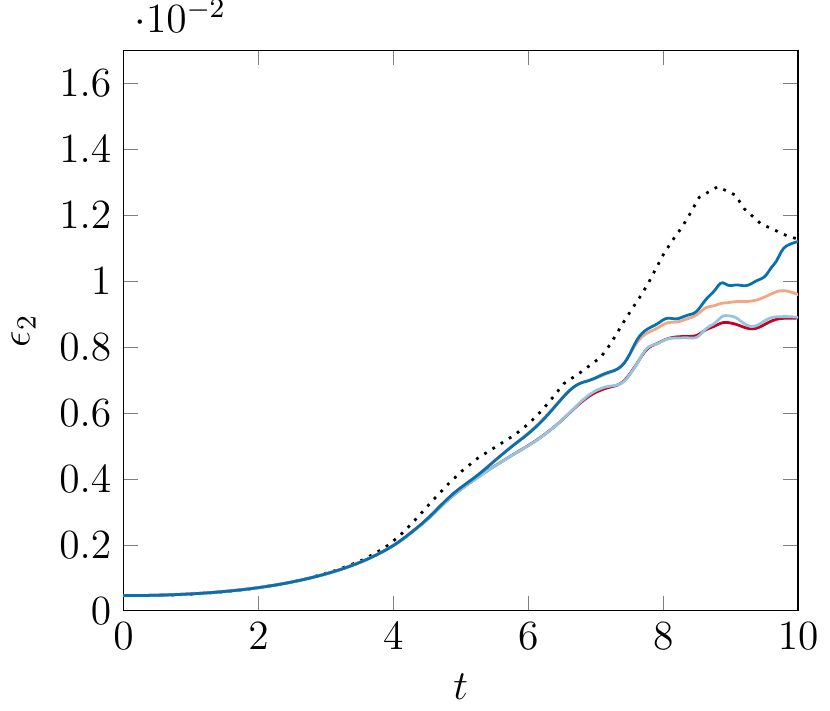}
			\caption{$M_a=0.31$}
			\label{fig:FRp4_fp_1600_m31}
		\end{subfigure}
		~
		\begin{subfigure}[b]{0.5\linewidth}
			\centering
			\includegraphics[width=\linewidth]{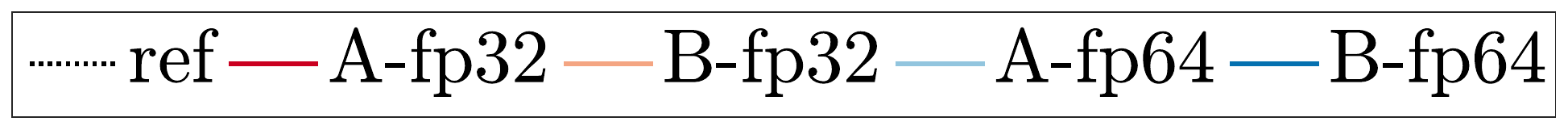}
		\end{subfigure}
		\caption[Taylor-Green Vortex Enstophy $R_e=1600$, $p=4$ and $80^3$ DoF for different variable storage precisions.]{Enstrophy of the Taylor-Green Vortex with $R_e=1600$, $p=4$ and $80^3$ degrees of freedom for storage methods A and B in 32 (fp32) and 64 (fp64) bit precision.}
		\label{fig:FR_fp_1600}
	\end{figure}	
	
	Finally, we investigate the numerical impact of varying the working precision of the calculation as applied to turbulent and transitional flows. For this, we limit our comparison to methods, A and B, as it was previously been shown that the largest differences were between these two methods.  The results of tests are shown in Fig.~\ref{fig:FR_fp_1600}, where 32-bit floating point (fp32) and 64-bit floating point (fp64) precisions were used. It is clear that the largest impact of changing the precision is at low Mach number. Coupled to the larger difference at higher Mach number being due to the variables stored, we believe that this is showing that at lower Mach number the scheme is more sensitive to numerical aliasing occurring in the interpolation. As the Mach number is increased and the physics begins to exhibit non-constant $\rho$, the small floating point errors in variables is more compatible with the physics and hence its effect appears to be lessened. This investigation into precision is of course limited, with the steeper gradients of discontinuities or solid boundaries likely to increase the impact.

%% file: nomenclature.tex
%!TEX root = ./aliasing_main.tex
\section{Nomenclature}
		\begin{tabbing}
			XXXXXXXX \= \kill% this line sets tab stop  
			\textit{Roman}\\
			A \> scheme storing $\mathbf{Q}_p$ and $\nabla\mathbf{Q}_p$ \\
			B \> scheme storing $\mathbf{Q}_c$ and $\nabla\mathbf{Q}_p$ \\
			C \> scheme storing $\mathbf{Q}_c$ and $\nabla\mathbf{Q}_c$ \\
			D \> scheme storing $\mathbf{Q}_{c+p}$ and $\nabla\mathbf{Q}_p$ \\
			$e_a$ \> aliasing error term \\
			$E_k$ \> volume averaged kinetic energy \\
			$\mathbf{f}^\mathrm{inv},\dots$ \> inviscid flux vector in $x,\dots$ \\
			$\mathbf{f}^\mathrm{vis},\dots$ \> viscous flux vector in $x,\dots$\\
			$h_L$~\&~$h_R$ \> left and right correction function \\
			$l_i$ \> $i^\mathrm{th}$ Lagrange basis polynomial \\ 
			$M_a$ \> Mach number \\
			$P_r$ \> Prandlt number \\
			$q_{p}$ \> non-normalised $p^\mathrm{th}$ order Lagrange basis \\
			$q_{p}^i$ \> non-normalised $p-1^\mathrm{th}$ order Lagrange basis formed from $q_p$ excluding $i^\mathrm{th}$ term\\
			$\mathbf{Q}_c$ \> conserved variables \\
			$\mathbf{Q}_{c+p}$ \> conserved variables with $E$ exchanged for $p$ \\
			$\mathbf{Q}_p$ \> primitive variables \\
			$\nabla\mathbf{Q}_c$ \> gradient of conserved variables \\
			$\nabla\mathbf{Q}_p$ \> gradient of primitive variables \\
			$R_e$ \> Reynolds number \\
			$T$ \> temperature \\
			$\nabla T$ \> gradient of temperature \\
			$\mathbf{V}$ \> vector of velocity components \\
			\\
		  	\textit{Greek}\\
		  	$\beta$ \> Icentropic Convecting Vortex spread rate \\
		  	$\Gamma_n(x)$ \> projection operator from real to reference domain, $\Gamma_n: \mathbf{\Omega}_n \mapsto \hb{\Omega}$ \\
		  	$\epsilon_1$ \> global averaged kinetic energy based dissipation, $-\mathrm{d}E_k/\mathrm{d}t$ \\ 
		  	$\epsilon_2$ \> global enstrophy based dissipation \\
		  	$\mu$ \> dynamic viscosity \\
		  	$\xi$ \> spatial variable in reference domain \\
		  	$\tau_{xx},\dots$ \> viscous stress tensor \\
		  	$\psi_i$ \> $i^\mathrm{th}$ order Legendre polynomial of the first kind \\
		  	$\pmb{\omega}$ \> vorticity, $\nabla\times\mathbf{V}$ \\
		  	$\mathbf{\Omega}$ \> spatial domain \\
		  	$\mathbf{\Omega}_n$ \> $n^\mathrm{th}$ spatial sub-domain \\
		  	$\hb{\Omega}$ \> reference domain \\
			\\
			\textit{Superscript} \\
			$\bullet^{\delta}$ \> approximation of variable in sub-domain\\
			$\bullet^{\delta C}$ \> correction to approximation of variable in sub-domain\\
			$\bullet^{\delta D}$ \> discontinuous approximation of variable in sub-domain\\
			$\bullet^{\delta I}$ \> common interface values based on approximation of variable in sub-domain\\
			$\hat{\bullet}$ \> variable transformed into reference domain\\
			$\tilde{\bullet}$ \> variable transformed into polynomial space\\
			$\bullet^{(n)}$ \> $n^\mathrm{th}$ derivative of variable \\
			\\
			\textit{Subscript} \\
			$\bullet_L$ \> variable at left interface \\
			$\bullet_R$ \> variable at right interface \\
			$\bullet_{x}$ \> differentiation of variable with respect to $x$ \\
			\\
			\textit{Other Symbols} \\
			$\mathcal{L}_n$ \> $n^\mathrm{th}$ order interpolation operator \\
			$\mathbb{N}$ \> set of natural numbers, i.e positive non-zero integers \\
			$\mathcal{O}$ \> big O notation of leading order in limiting behaviour \\ 
			$\mathbb{R}$ \> set of real numbers \\
			$\mathcal{R}_n$ \> $n^\mathrm{th}$ order interpolation remainder operator, i.e. $f-\mathcal{L}_nf$ \\

		\end{tabbing}